\documentclass[letterpaper, 10 pt, conference]{ieeeconf}  

\IEEEoverridecommandlockouts                              

\usepackage{cite}
\usepackage{amsmath,amssymb,amsfonts}
\usepackage[short]{optidef}
\usepackage{algorithmic}
\usepackage{graphicx}
\usepackage{textcomp}
\usepackage{xcolor}

\usepackage{pgfplots}
\usetikzlibrary{calc}
\usepackage{url}

\usepackage{hyperref}
\hypersetup{
	colorlinks,
	linkcolor={red!65!black},
	citecolor={green!65!black},
	urlcolor={blue!65!black}
}
\usepackage{tabularx}
\usepackage{hhline}
\usepackage{subcaption}
\newcommand{\bs}[1]{\boldsymbol{#1}}
\newcommand{\btx}[2]{\bs{#1}_{\mathrm{#2}}}
\newcommand{\qtx}[1]{\btx{q}{#1}}
\newcommand{\xtx}[1]{\btx{x}{#1}}

\newcommand{\bno}[2]{\bs{#1}_{#2}}
\newcommand{\zno}[1]{\bno{z}{#1}}
\newcommand{\uno}[1]{\bno{u}{#1}}

\newcommand{\hT}{^\mathrm{T}}

\newcommand{\Dt}{\Delta t}

\newcommand\ddfrac[2]{\frac{\displaystyle #1}{\displaystyle #2}}
\newcommand{\derivb}[2]{\ddfrac{\partial#1}{\partial#2}}

\title{\LARGE \bf
Linear-Quadratic Optimal Control in Maximal Coordinates
}

\author{Jan Br\"udigam$^{1}$ and Zachary Manchester$^{2}$
\thanks{$^{1}$Jan Br\"udigam is with the Department of Electrical and Computer Engineering, Technical University of Munich, Munich, Germany
        {\tt\small jan.bruedigam@tum.de}}%
\thanks{$^{2}$Zachary Manchester is with The Robotics Institute, Carnegie Mellon University, Pittsburgh, USA
        {\tt\small zacm@cmu.edu}}%
}

\begin{document}

\maketitle
\thispagestyle{empty}
\pagestyle{empty}

\begin{abstract}
The linear-quadratic regulator (LQR) is an efficient control method for linear and linearized systems. Typically, LQR is implemented in minimal coordinates (also called generalized or ``joint'' coordinates). However, other coordinates are possible and recent research suggests that there may be numerical and control-theoretic advantages when using higher-dimensional non-minimal state parameterizations for dynamical systems. One such parameterization is maximal coordinates, in which each link in a multi-body system is parameterized by its full six degrees of freedom and joints between links are modeled with algebraic constraints. Such constraints can also represent closed kinematic loops or contact with the environment. This paper investigates the difference between minimal- and maximal-coordinate LQR control laws. A case study of applying LQR to a simple pendulum and simulations comparing the basins of attraction and tracking performance of minimal- and maximal-coordinate LQR controllers suggest that maximal-coordinate LQR achieves greater robustness and improved tracking performance compared to minimal-coordinate LQR when applied to nonlinear systems.
\end{abstract}

\section{Introduction}\label{sec:intro}
Minimal coordinates (also called generalized or ``joint'' coordinates) have historically dominated robotic simulation and control, possibly due to the perception that they lead to greater computational efficiency. However, rigid body dynamics in maximal coordinates with Lagrange multipliers can be computed with similar efficiency as unconstrained dynamics in minimal coordinates \cite{baraff_linear-time_1996, brudigam_linear-time_2020}, and a substantial body of recent work, for example \cite{abraham_model-based_2017,suh_surprising_2020,bevanda_koopman_2021}, has suggested that higher-dimensional linear models may have more descriptive power than minimal models when approximating nonlinear systems.

Maximal-coordinate models of robotic systems are mathematically expressed as differential-algebraic equations (DAEs) in continuous time or algebraic difference equations (ADEs) in discrete time. Unfortunately, the classical derivation of the linear-quadratic regulator (LQR) is not well suited to such systems: If applied naively without explicitly accounting for ``joint'' or ``manifold'' constraints in the dynamics, the maximal-coordinate system will typically be mathematically uncontrollable, even if an equivalent minimal-coordinate realization is controllable. It is important to emphasize that the constraints in maximal-coordinate systems are not artificial state or input constraints that can be controlled by an actuator. Rather, they are hard mechanical constraints that are enforced by the physical structure. Such constraints always arise when using maximal coordinates, but they can also arise in minimal-coordinate settings, for example, for structures with closed-kinematic loops, during contact with the environment, or in nonholonomic settings.

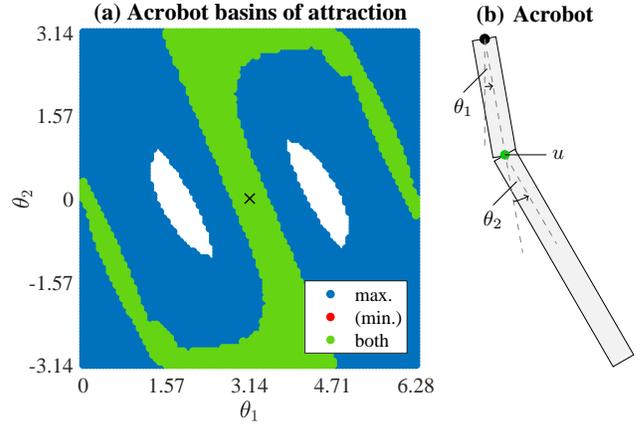
\begin{figure}[t] 
	\centering
	\resizebox{0.32\textwidth}{!}{\includegraphics{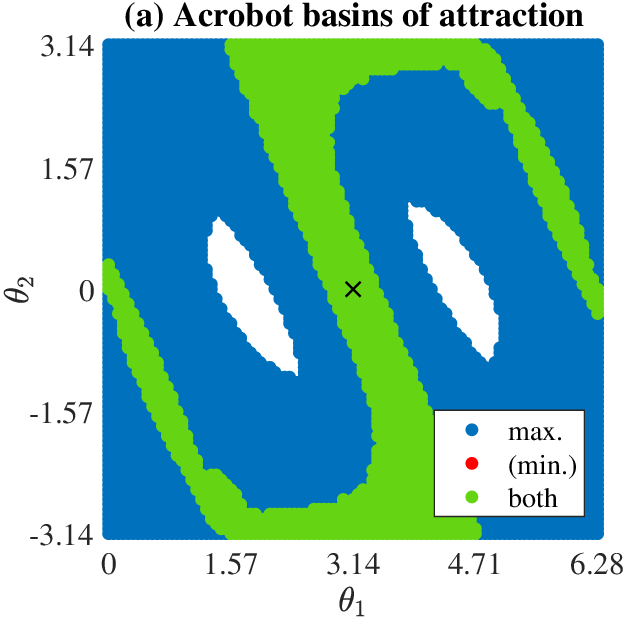}}
	\resizebox{0.14\textwidth}{!}{\begin{tikzpicture}
		\coordinate (Lw) at (0.0,0.0);
		\draw[fill=black!5,rotate=100] ($(Lw)+(0,0.0)$) rectangle ++(2,-0.4);
		\draw[fill=black!5,rotate=120] ($(Lw)+(-4.06,0.01)$) rectangle ++(4,-0.4);
		\draw[fill=black] ($(Lw)+(-0.14,2.0)$) circle (0.08);
		\fill[green!80!black] ($(Lw)+(0.2,0.05)$) circle (0.08);
		
		\draw[dashed,draw=black!50] ($(Lw)+(-0.14,2.0)$) -- ($(Lw)+(-0.14,0.2)$);
		\draw[dashed,draw=black!50] ($(Lw)+(-0.14,2.0)$) -- ($(Lw)+(0.5,-1.6)$);
		\draw[dashed,draw=black!50] ($(Lw)+(0.2,0.05)$) -- ($(Lw)+(1.07,-1.45)$);
		
		\draw[->] ($(Lw)+(-0.14,2.0)+(0.0,-0.8)$) arc (270:280:0.8);
		\draw[->] ($(Lw)+(0.2,0.05)+(0.8*0.174,-0.8*0.985)$) arc (280:300:0.8);
		
		\draw[-,draw=black] ($(Lw)+(-0.1,1.5)$) -- ($(Lw)+(-0.5,1.1)$);
		\draw[-,draw=black] ($(Lw)+(0.2,0.05)$) -- ($(Lw)+(0.9,0.05)$);
		\draw[-,draw=black] ($(Lw)+(0.4,-0.4)$) -- ($(Lw)+(0.0,-0.8)$);
		
		\node at ($(Lw)+(-0.5,0.85)$) {$\theta_1$};
		\node at ($(Lw)+(1.1,0.05)$) {$u$};
		\node at ($(Lw)+(0.0,-1.05)$) {$\theta_2$};
		
		\fill[white] ($(Lw)+(0.0,-4.5)$) circle (0.08); 
		
		\node[font=\bfseries,scale=1.1] at ($(Lw)+(0.7,2.4)$) {(b) Acrobot};
		\end{tikzpicture}}
	\caption{Basins of attraction for maximal- and minimal-coordinate LQR when stabilizing the upright equilibrium (black cross) of an acrobot from different initial configurations. The maximal-coordinate basin is shown in blue, the minimal-coordinate basin in red, and the basin for both in green. For this example, the maximal-coordinate basin includes the entire minimal-coordinate basin and is significantly larger.}
	\label{fig:resintro}
\end{figure}

In general, there are several approaches to control such mechanically constrained systems. Manipulator contact has been treated by constraining the end effector to always be in contact with the environment and then decomposing the system into subsystems relating to motion on and off the constraint manifold \cite{mills_force_1989,khayati_multistage_2006}. Contact arising with walking robots has been handled by deploying specialized control algorithms on the constraint manifold \cite{townsend_optimal_1977,posa_optimization_2016}. Systems with more general holonomic and nonholonomic constraints have been treated in a similar manner by either splitting the system into subsystems \cite{you_tracking_1993} or by analyzing the behavior on the constraint manifold \cite{guanfeng_liu_unified_2002}. Generally, the described methods are based on minimal (or reduced) coordinate representations and aim at reducing, splitting up, or eliminating the effect of constraints. A more mathematical perspective on mechanically constrained systems can be obtained by treating them as differential-algebraic equations (DAEs). General control of DAEs can be achieved with a variety of nonlinear feedback laws \cite{kumar_feedback_1996,krishnan_tracking_1994}. Optimal control of nonlinear DAEs has also been proposed \cite{kunkel_optimal_2008,gerdts_direct_2003}, and methods for systems described by linear DAEs have been developed as well \cite{bender_linear-quadratic_1987,kunkel_linear_1997}. For the purpose of deploying linear control methods, mechanically constrained systems and DAEs can also be linearized \cite{marz_linear_1995}. One common method for the linearization of mechanical systems described in minimal coordinates is to split the coordinates into independent and dependent variables \cite{kang_force_2003,negrut_practical_2006}. For an overview of more general linearization schemes for arbitrary coordinates in continuous time see \cite{gonzalez_assessment_2017}. 

The linear-quadratic regulator is a well-established control scheme for (locally) linear systems. However, LQR variants are typically used with minimal-coordinate representations and mechanical constraints are eliminated before deriving a controller \cite{mason_balancing_2016,savin_modification_2017,johnson_structured_2015}. There is also a wide variety of state- and input-constrained LQR control laws, for example \cite{scokaert_constrained_1998,mare_solution_2007,nguyen_fast_2016,ferranti_constrained_2016,laine_efficient_2019}, but these constraints are ``virtual.'' As such, they can be physically violated, and part of the control task is to ensure constraint satisfaction by an appropriate choice of inputs. As stated above, the nature of mechanical constraints and maximal-coordinate systems treated in this paper is fundamentally different.

Maximal coordinates represent the six degrees of freedom of each body in an articulated structure and the connections between bodies are expressed by explicit joint constraints. Several methods have been developed to efficiently compute dynamics of systems described in maximal coordinates, for example \cite{baraff_linear-time_1996,brudigam_linear-time_2020}. A related approach for trajectory optimization can be found in \cite{knemeyer_minor_2020}. Control of systems described in maximal coordinates has been mainly developed for individual vehicles, such as underwater, land, or aerial vehicles \cite{yilin_zhao_kinematics_1992,mellinger_minimum_2011,peng_constrained_2019}. In contrast, control of articulated systems parameterized in maximal coordinates has found very little attention. Some ideas related to this concept can, however, be found in cooperative robotics, where several individual vehicles form a virtual articulated structure when performing a common task \cite{lewis_high_1997,deshpande_distributed_2011}.

Currently, there exist some derivations of continuous-time LQR for DAEs \cite{bender_linear-quadratic_1987,kunkel_linear_1997}, and we will add to these concepts by providing a discrete-time LQR derivation for ADEs in Sec. \ref{sec:background_riccati_max}. The resulting control law will subsequently be used in Sec. \ref{sec:casestudy} for the case study of a simple pendulum to demonstrate that a linear maximal-coordinate control law can be interpreted as a nonlinear control law in minimal coordinates leading to potentially improved performance. This insight is further supported in Sec. \ref{sec:simulation} by comparing the basins of attraction of time-invariant LQR feedback laws in minimal and maximal coordinates, see Fig. \ref{fig:resintro} for an example, and by a comparison of trajectory tracking performance with time-varying LQR in minimal and maximal coordinates. Lastly, Sec. \ref{sec:discussion} provides a discussion and Sec. \ref{sec:conclusions} summarizes our conclusions.

\section{Background}\label{sec:background}
We will briefly depict the treatment of rigid body dynamics in maximal coordinates (see \cite{baraff_linear-time_1996} or \cite{brudigam_linear-time_2020} for details). Subsequently, we review the derivation of the classical unconstrained linear-quadratic regulator (LQR) from the dynamic programming principle with Riccati recursion (see  \cite{bertsekas_dynamic_1995} for details), and then state the Riccati recursion for LQR in maximal coordinates.

\subsection{Maximal Coordinates}
A single rigid body in space can be described by a position $\bs{x} \in \mathbb{R}^3$, a velocity $\bs{v} \in \mathbb{R}^3$, an orientation (unit quaternion) $\bs{q} \in \mathbb{R}^4$, and an angular velocity $\bs{\omega} \in \mathbb{R}^3$. We can group these quantities into a vector
\begin{equation}
\bs{z} = \begin{bmatrix}
\bs{x}\hT & \bs{v}\hT & \bs{q}\hT & \bs{\omega}\hT
\end{bmatrix}\hT.
\end{equation}

\begin{figure}
	\centering
	\begin{tikzpicture}
	\draw[fill=black!5,rounded corners=6pt,rotate=20]
	(0,0.0) rectangle ++(2,-0.4);
	\draw[fill=black!5,rounded corners=6pt,rotate=-30]
	(1.12,1.45) rectangle ++(2,-0.4);
	
	\draw[fill=black] (1.76,0.43) circle (0.065);
	
	\coordinate (l1) at (0.85,0.1);
	\draw[->,draw=blue!80,thick] (l1) -- ($(l1)+(-0.16,0.4)$);
	\draw[->,draw=blue!80,thick] (l1) -- ($(l1)+(0.4,0.16)$);
	\draw[->,draw=blue!80,thick] (l1) -- ($(l1)+(-0.15,-0.3)$);
	
	\coordinate (l2) at (2.47,-0.02);
	\draw[->,draw=blue!80,thick] (l2) -- ($(l2)+(0.23,0.36)$);
	\draw[->,draw=blue!80,thick] (l2) -- ($(l2)+(0.36,-0.23)$);
	\draw[->,draw=blue!80,thick] (l2) -- ($(l2)+(-0.15,-0.3)$);
	
	
	\coordinate (w) at (1.5,-1.2);
	\draw[->,draw=blue!80,thick] (w) -- ($(w)+(0.0,0.42)$);
	\draw[->,draw=blue!80,thick] (w) -- ($(w)+(0.42,0.0)$);
	\draw[->,draw=blue!80,thick] (w) -- ($(w)+(-0.15,-0.3)$);
	
	
	\node at (-0.2,0.5) {Link a};
	\node at (3.7,0.2) {Link b};
	
	
	\node at ($(w)+(1.0,-0.3)$) {Global frame};
	
	\draw ($(w)+(-0.07,0.07)$) to[out=150,in=-80 ] ($(l1)+(0.05,-0.05)$);
	\draw ($(w)+(0.07,0.07)$) to[out=50,in=180 ] ($(l2)+(-0.07,-0.02)$);
	
	\node at (0.5,-0.8) {$\xtx{a}$, $\qtx{a}$};
	\node at (2.5,-0.7) {$\xtx{b}$, $\qtx{b}$};
	\end{tikzpicture}
	\caption{Two links connected by a joint. Adopted from \cite{brudigam_linear-time_2020}.}
	\label{fig:twolinks}
\end{figure}
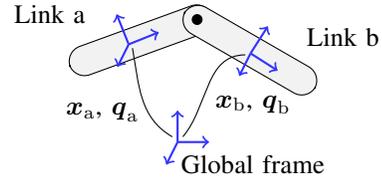

One or more bodies can be subject to constraints $\bs{g}$ which can, for example, represent joints. One example of such a physical constraint is the revolute joint connecting links a and b 
in Fig. \ref{fig:twolinks}. Such constraints on two bodies can generally be written as implicit equations
\begin{equation}
\bs{g}(\bs{z}_{\mathrm{a}},\bs{z}_{\mathrm{b}}) = \bs{0}.
\end{equation}

For each body, we can generally write the unconstrained dynamics as implicit equations in discrete time:
\begin{equation}
\bs{d}_0(\zno{k},\zno{k+1},\uno{k}) = \bs{0},
\end{equation}
where $k$ is the current time step and $\uno{k}$ are control inputs to the system, assuming a zero-order hold on the controls.

Constraints can be treated by introducing Lagrange multipliers (constraint forces) $\bs{\lambda}$ into the dynamics equations:
\begin{equation}\label{eq:dynamics}
\bs{d}(\zno{k},\zno{k+1},\uno{k},\bs{\lambda}_k) = \bs{d}_0 - G\hT \bs{\lambda}_k = \bs{0},\\
\end{equation}
where $G$ is the Jacobian of constraints $\bs{g}$ with respect to $\bs{z}$.

The dynamics equations $\bs{d}$ and the constraints $\bs{g}$ of multiple bodies form a system of algebraic difference equations (ADEs).

\subsection{Classical Riccati Recursion}\label{sec:background_riccati}

The linear-quadratic regulator (LQR) can be derived for four different scenarios: finite or infinite horizon, and discrete or continuous time. The Riccati recursion is one method to derive LQR for all four scenarios, and, while not always the most efficient method for each scenario, we will build on it for the remainder of the paper for clarity.

For discrete-time systems with a finite horizon, the Riccati recursion produces feedback gains at each time step. Continuing the recursion further leads to convergence of the feedback gains, which is simply the solution to the infinite-horizon problem. Furthermore, by taking the limit of the discrete time step as it goes to zero, we can recover the continuous time feedback gains for both finite-horizon and infinite-horizon settings.

We will provide the basic derivation of the classical recursion here to illustrate the general procedure.

The discrete-time optimal-control problem for LQR with states $\zno{k}$, controls $\uno{k}$, and weight matrices $Q$ and $R$ is:
\begin{mini}[2]
	{}{J_0 :=  \frac{1}{2}\zno{N}\hT Q \zno{N} + \frac{1}{2}\sum_{k=0}^{N-1}\left(\zno{k}\hT Q \zno{k} +  \uno{k}\hT R \uno{k} \right)}{}{}
	\addConstraint{\zno{k+1} = A\zno{k} + B\uno{k},}
	\label{eqn:optprob}
\end{mini}
where the equation $\zno{k+1} = A\zno{k} + B\uno{k}$ describes the unconstrained linear dynamics.

The cost-to-go function at time step $k$ is then given by
\begin{equation}\label{eqn:costtogo}
J_k = \frac{1}{2}\left(\zno{k}\hT Q \zno{k} +  \uno{k}\hT R \uno{k} \right) + J^{\star}_{k+1},
\end{equation}
with the optimal cost-to-go $J_k^{\star} = \min J_k$.

A state feedback law $\bs{u}$ depending linearly on the state $\bs{z}$ requires $J_k^{\star}$ to be of the form
\begin{equation}\label{eqn:Jkopt}
J_k^{\star} = \frac{1}{2}\zno{k}\hT P_k \zno{k}.
\end{equation}
Inserting \eqref{eqn:Jkopt} into \eqref{eqn:costtogo} yields
\begin{equation}\label{eqn:costtogoexpl}
\begin{alignedat}{1}
J_k &= \frac{1}{2}\Bigl(\zno{k}\hT  Q \zno{k} +  \uno{k}\hT R \uno{k}\Bigr. \\
&~~~~~~~~~~~ + \Bigl.\left(A\zno{k} + B\uno{k}\right)\hT P_{k+1} \left(A\zno{k} + B\uno{k}\right)\Bigr).
\end{alignedat}
\end{equation}

The optimal control $\uno{k}^{\star}$ can be found by taking the derivative of $J_k$ with respect to $\uno{k}$ and setting it equal to zero:
\begin{equation}\label{eqn:Jgrad}
\nabla_{\uno{k}} J_k = R\uno{k} + B\hT P_{k+1} \left(A\zno{k} + B\uno{k}\right) = 0.
\end{equation}
Rearranging yields
\begin{subequations}\label{eqn:uopt}
	\begin{align}
	\uno{k}^{\star} &= -\left(R + B\hT P_{k+1} B\right)^{-1}B\hT P_{k+1}A\zno{k}\\
	&= -K_k \zno{k}.
	\end{align}
\end{subequations}
Setting \eqref{eqn:Jkopt} equal to \eqref{eqn:costtogoexpl} after inserting \eqref{eqn:uopt} yields an update rule for $P_k$ after some rearranging:
\begin{equation}\label{eqn:Prec}
P_k = Q + K_k\hT R K_k + \bar{A}_{k}\hT P_{k+1}\bar{A}_{k},
\end{equation}
where $P_N = Q$ and
\begin{equation}
\bar{A}_{k} = A-BK_k.
\end{equation}

\subsection{Maximal-Coordinate Riccati Recursion}\label{sec:background_riccati_max}
We now provide the derivation of the linear-quadratic regulator (LQR) for maximal-coordinate systems and general ADE systems.

Given an implicit discrete-time dynamical system with implicit constraints
\begin{subequations}\label{eqn:nonlinear_dynamics}
	\begin{align}
	\bs{d}(\zno{k},\zno{k+1},\uno{k},\bs{\lambda}_k) = \bs{0},\\
	\bs{g}(\zno{k+1}) = \bs{0},
	\end{align}
\end{subequations}
we can apply the implicit function theorem to obtain linearized dynamics and constraints
\begin{subequations}\label{eqn:linear_dynamics}
	\begin{align}
	\zno{k+1} &= A\zno{k} + B\uno{k} + C\bs{\lambda}_k,\\
	G\zno{k+1} &= \bs{0}.
	\end{align}
\end{subequations}

Assuming linearized discrete-time dynamics and constraints of the form \eqref{eqn:linear_dynamics}, we can extend the classical Ricatti recursion. The optimal control problem (cf. \eqref{eqn:optprob}) now includes additional constraints:
\begin{mini}[2]
	{}{J_0 :=  \frac{1}{2}\zno{N}\hT Q \zno{N} + \frac{1}{2}\sum_{k=0}^{N-1}\left(\zno{k}\hT Q \zno{k} +  \uno{k}\hT R \uno{k} \right)}{}{}
	\addConstraint{\zno{k+1} = A\zno{k} + B\uno{k} + C\bs{\lambda}_k}
	\addConstraint{G\zno{k+1}=G\left(A\zno{k} + B\uno{k} + C\bs{\lambda}_k\right) = \bs{0}.}
	\label{eqn:optprobcon}
\end{mini}

The structures of the cost-to-go function \eqref{eqn:costtogo} and the optimal cost-to-go \eqref{eqn:Jkopt} are unchanged. However, inserting \eqref{eqn:Jkopt} into \eqref{eqn:costtogo} (cf. \eqref{eqn:costtogoexpl}) now includes the constrained update rule for $\zno{k+1}$:
\begin{subequations}\label{eqn:costtogoexplcon}
	\begin{align}
	\begin{split}
	&J_k = \frac{1}{2}\Bigl(\zno{k}\hT Q \zno{k} +  \uno{k}\hT R \uno{k}\Bigr.\\
	&~~ + \Bigl.\left(A\zno{k} + B\uno{k} + C\bs{\lambda}_k\right)\hT P_{k+1} \left(A\zno{k} + B\uno{k} + C\bs{\lambda}_k\right)\Bigr),
	\end{split}\label{eqn:costtogoexplconJ}\\
	&G\left(A\zno{k} + B\uno{k} + C\bs{\lambda}_k\right) = \bs{0}.\label{eqn:implcon}
	\end{align}
\end{subequations}

Instead of minimizing \eqref{eqn:costtogoexplcon} like a classical quadratic program with additional multipliers, we propose a more direct method to reduce the computational burden. Assuming linearly independent constraints---a common requirement for constrained optimization---the constraints \eqref{eqn:implcon} uniquely define the constraint forces $\bs{\lambda}_k(\uno{k})$ as implicit functions of $\uno{k}$ (and $\zno{k}$). Therefore, we obtain the necessary optimality condition by taking the gradient of \eqref{eqn:costtogoexplconJ} with respect to $\uno{k}$ (cf. \eqref{eqn:Jgrad}):
\begin{equation}
\nabla_{\uno{k}} J_k = R\uno{k} + D\hT P_{k+1}\left(A\zno{k} + B\uno{k} + C\bs{\lambda}_k\right) = \bs{0},
\end{equation}
where
\begin{equation}
D = \frac{\mathrm{d}\left(A\zno{k} + B\uno{k} + C\bs{\lambda}_k\right)}{\mathrm{d}\uno{k}} = B-C\left(GC\right)^{-1}GB,
\end{equation}
and $C${\scriptsize$\derivb{\bs{\lambda}_k}{\uno{k}}$} $= -C\left(GC\right)^{-1}GB$ is obtained by applying the implicit function theorem to the implicit constraints \eqref{eqn:implcon}.

For linearly independent constraints, we now have a non-singular linear system of equations at each time step:
\begin{equation}\label{eqn:linsyscon}
\begin{bmatrix}
R + D\hT P_{k+1} B & D\hT P_{k+1} C\\
GB & GC
\end{bmatrix}
\begin{bmatrix}
\uno{k}\\\bs{\lambda}_k
\end{bmatrix} = 
-\begin{bmatrix}
D\hT P_{k+1}\\G
\end{bmatrix}A\zno{k},
\end{equation}
from which we can obtain the optimal feedback laws for $\uno{k}$ and $\bs{\lambda}_k$ (cf. \eqref{eqn:uopt}) by inverting the matrix on the left hand side of \eqref{eqn:linsyscon}:
\begin{equation}\label{eqn:ulaopt}
\begin{bmatrix}
\uno{k}^{\star}\\\bs{\lambda}_k^{\star}
\end{bmatrix} = 
-\begin{bmatrix}
K_k\\L_k
\end{bmatrix}\zno{k}.
\end{equation}
The matrix $K$ contains the control gains, whereas $L$ models the mapping of constraint forces during the derivation. In the physical system, these forces are ``automatically'' applied, so, for the actual controller, only $K$ is retained.

As before, the recursive update rule for $P_k$ (cf. \eqref{eqn:Prec}) is
\begin{equation}\label{eqn:Preccon}
P_k = Q + K_k\hT R K_k + \bar{A}_k\hT P_{k+1}\bar{A}_k,
\end{equation}
where now
\begin{equation}
\bar{A}_k = A-BK_k-CL_k.
\end{equation}

As stated above, the infinite-horizon and continuous-time control laws can be recovered from this derivation.

\begin{figure*}[t] 
	\centering
	\begin{subfigure}{0.325\textwidth}
		\resizebox{\linewidth}{!}{\includegraphics{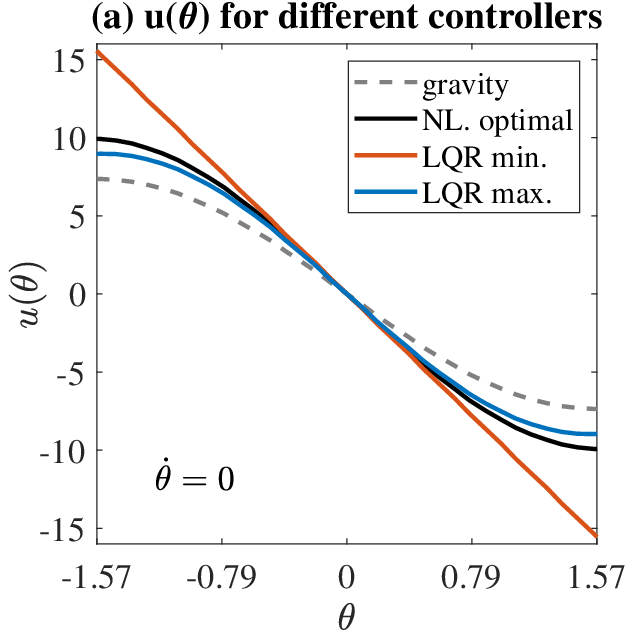}}
		\refstepcounter{subfigure}\label{fig:respend0}
	\end{subfigure}
	\begin{subfigure}{0.325\textwidth}
		\resizebox{\linewidth}{!}{\includegraphics{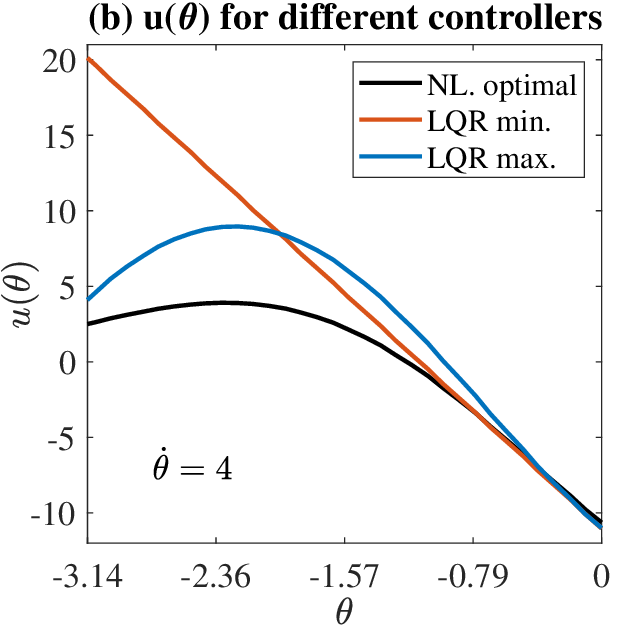}}
		\refstepcounter{subfigure}\label{fig:respend4}
	\end{subfigure}
	\begin{subfigure}{0.325\textwidth}
		\resizebox{\linewidth}{!}{\includegraphics{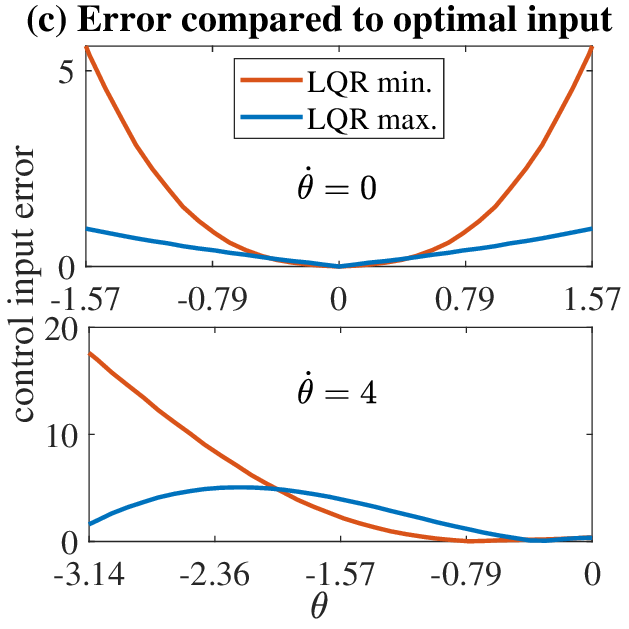}}
		\refstepcounter{subfigure}\label{fig:respenderror}
	\end{subfigure}
	\caption{Comparison of minimal- and maximal-coordinate LQR control laws with the optimal nonlinear control law. Maximal-coordinate control law (LQR max.) shown in blue, minimal-coordinate control law (LQR min.) in red, optimal control law (NL. optimal) in black. (a) Control laws for $\dot{\theta}=0$. (b) Control laws for $\dot{\theta}=4$. (c) Relative error of the LQR control laws compared to the optimal control law.}\label{fig:respend}
\end{figure*} 

\section{Case Study: Simple Pendulum}\label{sec:casestudy}

We demonstrate the difference between minimal- and maximal-coordinate LQR control laws by trying to stabilize a simple two-dimensional pendulum with mass $m=1$, length $l=1$, and inertia $J=\frac{1}{12}$ in the upright position. The minimal-coordinate states of the system are the angle $\theta$ and angular velocity $\dot{\theta}$. The maximal-coordinate states are the angle $\theta$, angular velocity $\dot{\theta}$, positions $x$ and $y$, and velocities $\dot{x}$ and $\dot{y}$ with gravity acting along the $\mathrm{y}$-axis.

In order to make the analysis of minimal- and maximal-coordinate LQR comparable, we only apply weights to minimal-coordinate states, even when using maximal coordinates. Additionally, the cost functions are defined such that, in a linearized sense, they have the same value in a region near the desired state. So, assuming a mapping $\bs{f}$ from maximal coordinates $\bs{z}$ to minimal coordinates $\bs{c}$
\begin{equation}
\bs{c} = \bs{f}(\bs{z}),
\end{equation}
the linearized relation between minimal- and maximal-coordinate cost at the reference point $\bs{z}^*$ is given as
\begin{equation}
\bs{c}\hT Q \bs{c} = \bs{z}\hT F\hT Q F \bs{z},
\end{equation}
where $F$ is the Jacobian of $\bs{f}$ and we choose $\bs{f}$ such that $\bs{f}(\bs{z}^*)=\bs{0}$. In the case of the pendulum, the cost functions are also globally exactly the same.

The desired upright state in minimal coordinates is $(\theta, \dot{\theta}) = (0,0)$, with nominal control input $u = 0$. As corresponding LQR weights we choose $Q_{\text{min}} = \mathrm{diag}(1,1)$ and $R = 1$. In maximal coordinates, the desired state is $(\theta, \dot{\theta}, x, \dot{x}, y, \dot{y}) = (0,0,0,0,0.5,0)$, and $Q_{\text{max}} = \mathrm{diag}(1,1,0,0,0,0)$, only applying weights to the minimal-coordinate states. The resulting minimal- and maximal-coordinate LQR feedback matrices for this setup are $K_{\text{min}} = [9.91, 2.76]$ and $K_{\text{max}} = [-19.30, 0, -4.13, 0, 0.44, 0.69]$, respectively. It is important to note that in maximal coordinates we have feedback terms for $x$ and $\dot{x}$ deviations because of the mechanical coupling despite not placing any weights on these states directly. On a real system, we would directly compute the actual position $x = -\frac{1}{2}\sin\left(\theta\right)$ and velocity $\dot{x} = -\frac{1}{2}\cos\left(\theta\right)\dot{\theta}$ from the measured angle $\theta$ and angular velocity $\dot{\theta}$ to calculate the control error. Therefore, if we express the maximal-coordinate control in terms of $\theta$ and $\dot{\theta}$, we can obtain the control law
\begin{equation*}
	u_{\text{max}}(\theta,\dot{\theta}) = -\frac{19.30}{2}\sin\left(\theta\right) -\frac{4.13}{2}\cos\left(\theta\right)\dot{\theta} - 0.44\theta - 0.69\dot{\theta},
\end{equation*}
which is nonlinear in $\theta$. Note that for small angles $\theta$, i.e. $\sin\theta\approx\theta$, $\cos\theta\approx1$, we get $u_{\text{max}}(\theta,\dot{\theta}) \approx -10.09\theta - 2.75\dot{\theta}$, which is the same as the minimal coordinate control law up to rounding errors.

In Fig. \ref{fig:respend} we compared the control laws with the optimal nonlinear control law for the same quadratic cost function. The nonlinear control law was obtained with the trajectory optimization tool \textit{ALTRO} \cite{howell_altro_2019}. Figure \ref{fig:respend0} shows the resulting control input for angles $\theta$ and a fixed angular velocity $\dot{\theta}=0$. We only show angles from $-\frac{\pi}{2}$ to $\frac{\pi}{2}$, since the optimal control law is discontinuous at slightly larger angles, i.e. the optimal path would be to swing the pendulum the other way around. As expected, for small angles, the nonlinear, minimal-coordinate, and maximal-coordinate control laws roughly match. However, for larger angles, the maximal-coordinate control law closely follows the nonlinear control law, while the minimal-coordinate control law deviates. We also plotted the gravitational torque at a given angle for reference, since all control laws have to drive the pendulum against gravity. Figure \ref{fig:respend4} shows the same analysis for angles from $-\pi$ to $0$ and a fixed angular velocity $\dot{\theta}=4$. As before, for small angles, the control laws roughly match, and again, for larger angles, the maximal-coordinate control law follows the nonlinear optimal control law more closely, even though initially it shows a slightly larger deviation than the minimal-coordinate control law. In Fig. \ref{fig:respenderror} we plotted the relative error compared to the nonlinear control law for the two scenarios for convenience.

\begin{figure*}[t]
	\centering
	\begin{subfigure}{0.3\textwidth}
		\begin{tikzpicture}
		\coordinate (Lw) at (0.0,0.0);
		\draw[fill=black!5,rotate=100] ($(Lw)+(0,0.0)$) rectangle ++(2,-0.4);
		\draw[fill=black!5,rotate=120] ($(Lw)+(-4.06,0.01)$) rectangle ++(4,-0.4);
		\draw[fill=black] ($(Lw)+(-0.14,2.0)$) circle (0.08);
		\fill[green!80!black] ($(Lw)+(0.2,0.05)$) circle (0.08);
		\draw[->] ($(Lw)+(0.2+0.4,0.05)$) arc (0:345:0.4);
		
		\coordinate (Cw) at ($(Lw)+(-0.14,2.0)$);
		\draw[->,draw=black!80,thick] (Cw) -- ($(Cw)+(0.0,0.5)$);
		\draw[->,draw=black!80,thick] (Cw) -- ($(Cw)+(0.5,0.0)$);
		\draw[->,draw=black!80,thick] (Cw) -- ($(Cw)+(-0.26,-0.35)$);
		\node at ($(Cw)+(-0.4,-0.2)$) {x};
		\node at ($(Cw)+(0.6,0.15)$) {y};
		\node at ($(Cw)+(0.2,0.4)$) {z};
		
		\draw[dashed,draw=black!50] (Cw) -- ($(Lw)+(-0.14,0.2)$);
		\draw[dashed,draw=black!50] (Cw) -- ($(Lw)+(0.5,-1.6)$);
		\draw[dashed,draw=black!50] ($(Lw)+(0.2,0.05)$) -- ($(Lw)+(1.07,-1.45)$);
		
		\draw[->] ($(Lw)+(-0.14,2.0)+(0.0,-0.8)$) arc (270:280:0.8);
		\draw[->] ($(Lw)+(0.2,0.05)+(0.8*0.174,-0.8*0.985)$) arc (280:300:0.8);
		
		\draw[-,draw=black] ($(Lw)+(-0.1,1.5)$) -- ($(Lw)+(-0.5,1.1)$);
		\draw[-,draw=black] ($(Lw)+(0.2,0.05)$) -- ($(Lw)+(0.9,0.3)$);
		\draw[-,draw=black] ($(Lw)+(0.4,-0.4)$) -- ($(Lw)+(0.0,-0.8)$);
		
		\node at ($(Lw)+(1.0,1.6)$) {Link 1};
		\node at ($(Lw)+(1.0,1.2)$) {$m_1,J_1,l_1$};
		\node at ($(Lw)+(2.4,-1.2)$) {Link 2};
		\node at ($(Lw)+(2.4,-1.6)$) {$m_2,J_2,l_2$};
		
		\node at ($(Lw)+(-0.5,0.85)$) {$\theta_1$};
		\node at ($(Lw)+(1.1,0.3)$) {$u$};
		\node at ($(Lw)+(0.0,-1.05)$) {$\theta_2$};
		
		\node[font=\bfseries] at ($(Lw)+(1.1,2.8)$) {(a) Acrobot};
		\end{tikzpicture}\refstepcounter{subfigure}\label{fig:acrobot}
	\end{subfigure}
	\begin{subfigure}{0.3\textwidth}
		\begin{tikzpicture}
		\coordinate (Lw) at (0.0,0.0);
		\coordinate (Cw) at ($(Lw)+(-1.0,-0.2)$);
		\fill[green!80!black] ($(Lw)+(-0.3,-0.15)$) rectangle ++(2.6,-0.1);
		\draw[->] ($(Lw)+(-0.3+2.6,-0.2)$) -- ($(Lw)+(-0.3+3.0,-0.2)$);
		\draw[fill=black!5,rotate=0] ($(Lw)+(0,0.0)$) rectangle ++(2,-0.4);
		\draw[fill=black!5,rotate=110] ($(Lw)+(-0.55,-0.67)$) rectangle ++(4,-0.4);
		\draw[fill=black] ($(Lw)+(1.0,-0.2)$) circle (0.08);

		\draw[->,draw=black!80,thick] (Cw) -- ($(Cw)+(0.0,0.5)$);
		\draw[->,draw=black!80,thick] (Cw) -- ($(Cw)+(0.5,0.0)$);
		\draw[->,draw=black!80,thick] (Cw) -- ($(Cw)+(-0.26,-0.35)$);
		\node at ($(Cw)+(-0.4,-0.2)$) {x};
		\node at ($(Cw)+(0.6,0.15)$) {y};
		\node at ($(Cw)+(0.2,0.4)$) {z};
		
		\draw[dashed,draw=black!50] ($(Lw)+(1.0,-0.8)$) -- ($(Lw)+(1.0,1.7)$);
		\draw[dashed,draw=black!50] (Cw) -- ($(Cw)+(0.0,-0.6)$);
		\draw[dashed,draw=black!50] ($(Lw)+(1.0,-0.2)$) -- ($(Lw)+(0.38,1.5)$);
		
		\draw[->] ($(Lw)+(1.0,0.0)+(0.0,0.8)$) arc (90:113:0.8);
		\draw[->] ($(Cw)+(0.0,-0.4)$) -- ($(Lw)+(1.0,-0.6)$);
		
		\draw[-,draw=black] ($(Lw)+(0.9,0.5)$) -- ($(Lw)+(1.5,0.8)$);
		\draw[-,draw=black] ($(Lw)+(2.15,-0.2)$) -- ($(Lw)+(2.35,-0.4)$);
		\draw[-,draw=black] ($(Lw)+(0.0,-0.6)$) -- ($(Lw)+(0.0,-0.8)$);
		
		\node at ($(Lw)+(1.0,3.4)$) {Pole};
		\node at ($(Lw)+(1.0,3.0)$) {$m_\mathrm{p},J_\mathrm{p},l_\mathrm{p}$};
		\node at ($(Lw)+(2.6,0.8)$) {Cart};
		\node at ($(Lw)+(2.6,0.4)$) {$m_\mathrm{c},l_\mathrm{c}$};
		
		\node at ($(Lw)+(1.65,0.95)$) {$\theta_\mathrm{p}$};
		\node at ($(Lw)+(2.35,-0.6)$) {$u$};
		\node at ($(Lw)+(0.0,-1.0)$) {$y_\mathrm{c}$};
		
		\fill[white] ($(Lw)+(0.0,4.75)$) circle (0.08); 
		
		\node[font=\bfseries] at ($(Lw)+(0.8,5.1)$) {(b) Cartpole};
		\end{tikzpicture}\refstepcounter{subfigure}\label{fig:cartpole}
	\end{subfigure}
	\begin{subfigure}{0.3\textwidth}
		\begin{tikzpicture}
		\coordinate (Lw) at (0.0,0.0);
		\coordinate (Cw) at ($(Lw)$);
		\node [shading = axis,rectangle, left color=darkgray!50!white, right color=white, shading angle=0, minimum width = 120] at ($(Cw)+(0.0,-0.13)$) {};
		\draw[fill=black!5,rotate=120.5] ($(Lw)+(0.0,0.2)$) rectangle ++(4,-0.4);
		\draw[fill=black!5,rotate=32] ($(Lw)+(0.0,0.2)$) rectangle ++(4,-0.4);
		\draw[fill=black!5,rotate=17] ($(Lw)+(-0.92,4.05)$) rectangle ++(2,-0.4);
		\draw[fill=black!5,rotate=110] ($(Lw)+(0.8,-3.7)$) rectangle ++(2,-0.4);
		\draw[fill=black!5,rotate=0] ($(Lw)+(-0.1,4.2)$) rectangle ++($(2*1.4142,-0.4)$);
		
		\draw[fill=black] ($(Lw)$) circle (0.08);
		\draw[fill=black] ($(Lw)+(-2.02,3.45)$) circle (0.08);
		\draw[fill=black] ($(Lw)+(3.4,2.1)$) circle (0.08);
		\fill[green!80!black] ($(Lw)+(-0.1,4.0)$) circle (0.08);
		\fill[green!80!black] ($(Lw)+(2*1.4142-0.1,4.0)$) circle (0.08);
		\draw[->] ($(Lw)+(-0.1+0.4,4.0)$) arc (0:345:0.4);
		\draw[->] ($(Lw)+(2*1.4142-0.1+0.4,4.0)$) arc (0:345:0.4);

		\draw[->,draw=black!80,thick] (Cw) -- ($(Cw)+(0.0,0.5)$);
		\draw[->,draw=black!80,thick] (Cw) -- ($(Cw)+(0.5,0.0)$);
		\draw[->,draw=black!80,thick] (Cw) -- ($(Cw)+(-0.26,-0.35)$);
		\node at ($(Cw)+(-0.4,-0.2)$) {x};
		\node at ($(Cw)+(0.55,0.25)$) {y};
		\node at ($(Cw)+(0.2,0.4)$) {z};

		\draw[dashed,draw=black!50] (Cw) -- ($(Cw)+(0.0,3.4)$);
		\draw[dashed,draw=black!50] (Cw) -- ($(Cw)+(0.9,0.0)$);
		\draw[dashed,draw=black!50] ($(Lw)+(1.4142-0.05,4.0)$) -- ($(Lw)+(1.4142-0.05,3.15)$);
		\draw[dashed,draw=black!50] ($(Lw)+(1.4142-0.05,4.0)$) -- ($(Lw)+(1.4142-0.05-0.8,4.0)$);
		\draw[->] ($(Cw)+(0.0,3.3)$) -- ($(Lw)+(1.4142-0.05,3.3)$);
		\draw[->] ($(Cw)+(0.8,0.0)$) -- ($(Lw)+(0.8,4.0)$);

		\node at ($(Lw)+(2.2,2.7)$) {Upper};
		\node at ($(Lw)+(2.2,2.3)$) {$m_\mathrm{u},J_\mathrm{u},l_\mathrm{u}$};
		\node at ($(Lw)+(2.6,0.9)$) {Lower};
		\node at ($(Lw)+(2.6,0.5)$) {$m_\mathrm{l},J_\mathrm{l},l_\mathrm{l}$};
		\node at ($(Lw)+(2.1,5.0)$) {Base};
		\node at ($(Lw)+(2.1,4.6)$) {$m_\mathrm{b},l_\mathrm{b}$};
		\draw[-,draw=black] ($(Lw)+(2.5,2.9)$) -- ($(Lw)+(2.9,3.5)$);
		\draw[-,draw=black] ($(Lw)+(1.7,2.7)$) -- ($(Lw)+(-1.4,3.5)$);
		\draw[-,draw=black] ($(Lw)+(2.6,1.1)$) -- ($(Lw)+(2.7,1.6)$);
		\draw[-,draw=black] ($(Lw)+(2.0,0.8)$) -- ($(Lw)+(-0.9,1.6)$);
		
		\draw[-,draw=black] ($(Cw)+(1.2,3.5)$) -- ($(Cw)+(1.2,3.3)$);
		\draw[-,draw=black] ($(Cw)+(0.6,3.62)$) -- ($(Cw)+(0.8,3.62)$);
		
		\draw[-,draw=black] ($(Lw)+(2*1.4142-0.1,4.0)$) -- ($(Lw)+(2*1.4142+0.5,4.0)$);
		\draw[-,draw=black] ($(Lw)+(-0.1,4.0)$) -- ($(Lw)+(-0.6,4.3)$);
		\node at ($(Lw)+(-0.85,4.4)$) {$u_1$};
		\node at ($(Lw)+(2*1.4142+0.75,4.0)$) {$u_2$};
		\node at ($(Cw)+(1.2,3.6)$) {$y_\mathrm{b}$};
		\node at ($(Cw)+(0.45,3.62)$) {$z_\mathrm{b}$};
		
		\node at ($(Cw)+(2.9,-0.1)$) {Ground};
		\draw[-,draw=black] ($(Cw)+(2.1,-0.1)$) -- ($(Cw)+(2.3,-0.1)$);
		
		\fill[white] ($(Lw)+(0.0,5.75)$) circle (0.08); 
		
		\node[font=\bfseries] at ($(Lw)+(0.6,5.9)$) {(c) Delta robot};
		\end{tikzpicture}\refstepcounter{subfigure}\label{fig:deltabot}
	\end{subfigure}
	\caption{Systems for the LQR basin of attraction analysis. Actuated joints highlighted in green. (a) Acrobot. (b) Cartpole. (c) Delta robot. The base angle of the delta robot is kept constant by parallel structures (not drawn).}
	\label{fig:systems}
\end{figure*}
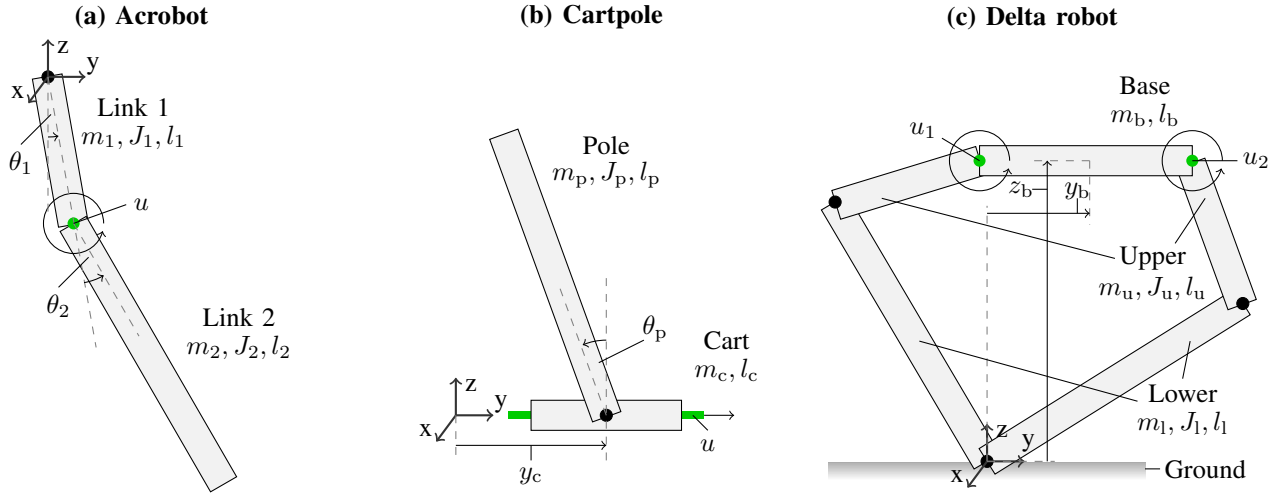

\section{Simulation}\label{sec:simulation}
The LQR feedback control was applied to more complex systems. As the dynamics simulation tool we chose the \textit{ConstrainedDynamics} package based on \cite{brudigam_linear-time_2020} written in the programming language \textit{Julia}. The code for our algorithm and all experiments is available at: \url{https://github.com/janbruedigam/ConstrainedControl.jl}.

The first part of the simulation analysis consists of determining the basins of attraction for minimal- and maximal-coordinate LQR for different systems. The infinite-horizon time-invariant linear controllers were calculated for a stationary reference state and closed-loop simulations were performed from densely sampled initial states. The stationary analysis is followed by a time-varying LQR tracking example. Given a nominal state and control trajectory, we calculated a linear time-varying model and the corresponding time-varying LQR. As in the pendulum case study, we only applied weights to minimal-coordinate states, even when using maximal coordinates, and we (linearly) matched the cost functions.

\subsection{Basin of Attraction Setup}
For the basin of attraction analysis of each system we chose initial conditions on a grid spanning the entire configuration space with zero initial velocities. We simulated each system for $25$ seconds with a time step of $\Dt = 1$ms. An initial condition was counted as inside the basin of attraction if the system satisfied $\lVert\bs{z}-\bs{z}^*\rVert<0.1$ within $25$ seconds. The control inputs were not limited. However, we stopped simulations if angular velocities exceeded $100\pi$ radians per second as the simulation no longer delivers reliable results for such excessive velocities. Velocities of this magnitude indicated instability.

For the analysis, we investigated the acrobot (Fig. \ref{fig:acrobot}), cartpole (Fig. \ref{fig:cartpole}), and a 2-dimensional delta robot with closed kinematic loops (Fig. \ref{fig:deltabot}). The control goal for all three systems was to stabilize an unstable equilibrium point at which the linearization was performed. All weight matrices $Q$ and $R$ are purely diagonal. The mechanical properties for all components are stated in Table \ref{tab:mech}.

\begin{table}
	\begin{center}
		\begin{tabularx}{\columnwidth}{l l l l l l l l l}
			\hline
			\hline
			& & Link 1 & Link 2 & Cart & Pole & Base & Upper & Lower \\
			\hline
			& $m$ & 1.0   & 1.0   & 0.5 & 1.0   & 0.71  & 0.5   & 1.0 \\
			& $J$ & 0.084 & 0.334 & -   & 0.084 & -     & 0.011 & 0.084 \\
			& $l$ & 1.0   & 2.0   & 0.5 & 1.0   & 0.71  & 0.5   & 1.0 \\
			\hline
			\hline
		\end{tabularx} 
		\caption{Mechanical parameters of the systems' components.}\label{tab:mech}	 
	\end{center}
\vspace*{-5mm}
\end{table}

\subsubsection{Acrobot Setup}
The acrobot is a double pendulum with a single actuator placed between the two links---mimicking the dynamics of an acrobat---resulting in an underactuated but controllable system. Generally, inverted double pendulums are very nonlinear systems which makes them interesting from a control perspective. The control goal for the acrobot is to stabilize it in the upright position.

The desired state in minimal coordinates is $(\theta_1, \theta_2, \dot{\theta}_1, \dot{\theta}_2) = (\pi,0,0,0)$ with nominal control input $u = 0$. The corresponding LQR weights are $Q = \mathrm{diag}(1,1,1,1)$ and $R = 1$.

\subsubsection{Cartpole Setup}
The cartpole consists of a pendulum (pole) mounted on a moving cart. The pendulum is connected by a passive joint. The actuated cart can only move along the $\mathrm{y}$-axis. The control goal for the cartpole is to stabilize the inverted pendulum while keeping the cart near the origin.

The desired state in minimal coordinates is $(y_\mathrm{c}, \theta_\mathrm{p}, \dot{y}_\mathrm{c}, \dot{\theta}_\mathrm{p}) = (0,0,0,0)$ with nominal control input $u = 0$. The corresponding LQR weights are $Q = \mathrm{diag}(1,1,1,1)$ and $R = 1$.

\begin{figure*}[t] 
	\centering
	\begin{subfigure}{0.325\textwidth}
		\resizebox{\linewidth}{!}{\includegraphics{plots/acrobot.eps}}
		\refstepcounter{subfigure}\label{fig:resacrobot}
	\end{subfigure}
	\begin{subfigure}{0.325\textwidth}
		\resizebox{\linewidth}{!}{\includegraphics{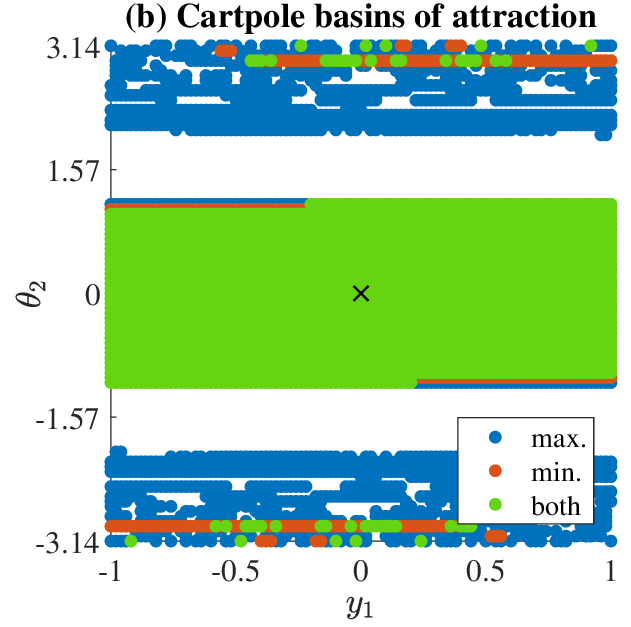}}
		\refstepcounter{subfigure}\label{fig:rescartpole}
	\end{subfigure}
	\begin{subfigure}{0.325\textwidth}
		\resizebox{\linewidth}{!}{\includegraphics{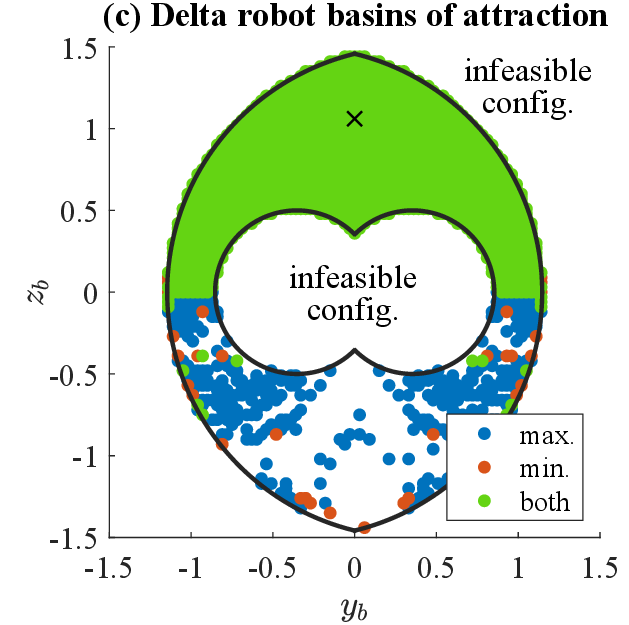}}
		\refstepcounter{subfigure}\label{fig:resdeltabot}
	\end{subfigure}
	\caption{Results for the LQR basin of attraction analysis. Maximal-coordinate basins (max.) shown in blue, minimal-coordinate basins (min.) in red, basins for both in green. The desired configurations are indicated by black crosses. (a) Results for acrobot. (b) Results for cartpole. (c) Results for the delta robot and indication of structurally infeasible initial configurations.}\label{fig:basins}
\end{figure*} 

\subsubsection{Delta Robot Setup}
Delta robots can be used to perform many tasks, including pick-and-place and milling. For the analysis we chose a 2-dimensional delta robot consisting of five links. The closed kinematic loops of the robot make a straight-forward application of minimal-coordinate LQR difficult. In contrast, the maximal-coordinate approach can handle this system without modifications to the algorithm. The system is driven by two actuators connected to the base resulting again in an underactuated but controllable system. For the control goal, we aimed at treating the delta robot as an inverted pendulum with the goal of stabilizing the system in an upright configuration close to the one depicted in Fig. \ref{fig:deltabot}. The system can also be interpreted as an instance of a bipedal robot with feet pinned to the ground, underactuated articulated legs (``Lower'' and ``Upper''), and a torso (``Base'') that is supposed to be balanced.

The desired state in minimal (reduced) coordinates is $(y_\mathrm{b}, z_\mathrm{b}, \dot{y}_\mathrm{b}, \dot{z}_\mathrm{b}) = (0,1.061,0,0)$ with nominal control inputs $(u_1,u_2) = (6.788,-6.788)$ compensating gravity. The corresponding LQR weights are $Q = \mathrm{diag}(100,100,1,1)$ and $R = 0.01$. For all initial conditions we set the ``legs'' of the system to be pointing outwards, as depicted in Fig. \ref{fig:deltabot}.

\subsection{Basin of Attraction Results}
The results for the simulations of the acrobot, cartpole, and delta robot are presented in Fig. \ref{fig:basins}. The regions for maximal coordinates are colored in blue, the regions for minimal coordinates in red. Regions for which both maximal- and minimal-coordinate LQR converged are colored in green.

\subsubsection{Acrobot Results}
The basin of attraction for the acrobot with minimal-coordinate LQR shows ``tilted'' regions around the fully upright ($\theta_1=\pi$, $\theta_2=0$) and fully hanging state ($\theta_1=0$, $\theta_2=0$) with connecting regions in between. A large area that is not stabilizable exists in between the tilted regions. In comparison, the stabilizable region for maximal-coordinate LQR has some resemblance in shape but is significantly larger. 

\subsubsection{Cartpole Results}
The cartpole's basin of attraction for minimal-coordinate LQR consists of a stabilizable region for initial pendulum angle deviations within $\theta_2 \approx \pm\frac{1}{3}\pi$. The position of the cart does not have a significant effect on the stability. A few additional stabilizable initial conditions for hanging pendulum configurations ($\theta_2\approx \pi$) were also detected. For angle deviations within $\theta_2\approx\pm\frac{1}{3}\pi$ maximal-coordinate LQR yields almost the same region of attraction. Additionally, some more stabilizable points can be found for configurations around the hanging state within $\theta_2 \approx \pi\pm\frac{1}{3}\pi$. 

\subsubsection{Delta Robot Results}
For z-positions above zero, both minimal- and maximal-coordinate LQR are capable of driving the delta robot to the desired state at any feasible initial configuration. The area below zero is less consistent, especially for minimal-coordinate LQR. In contrast, maximal-coordinate LQR is still able to drive the system into the desired position for several initial configurations.

\subsection{Tracking LQR Setup}			
Tracking of a nominal trajectory with time-varying LQR was performed on a triple-pendulum cartpole system. The mechanical properties of cart and poles are the same as for the single-pendulum cartpole (see Fig. \ref{fig:cartpole} and Table \ref{tab:mech}). The reference was a swing-up of the triple pendulum into the upright position, while LQR was deployed to correct for noise and model errors.

As before, the nominal state and control trajectory were calculated with \textit{ALTRO}. Subsequently, a time-varying model in minimal and maximal coordinates was calculated to obtain the respective LQR feedback laws. Again, the constant weight matrices $Q$ and $R$ only applied cost to minimal-coordinate states and the maximal-coordinate cost function was matched to first order. The tracked state in minimal coordinates is $(y_\mathrm{c}, \theta_\mathrm{p1}, \theta_\mathrm{p2}, \theta_\mathrm{p3}, \dot{y}_\mathrm{c}, \dot{\theta}_\mathrm{p1}, \dot{\theta}_\mathrm{p2}, \dot{\theta}_\mathrm{p3})$ and the cart's control input is $u$. The corresponding LQR weights are $Q = \mathrm{diag}(10,10,10,10,1,1,1,1)$ and $R = 0.1$.

\begin{figure*}[t]
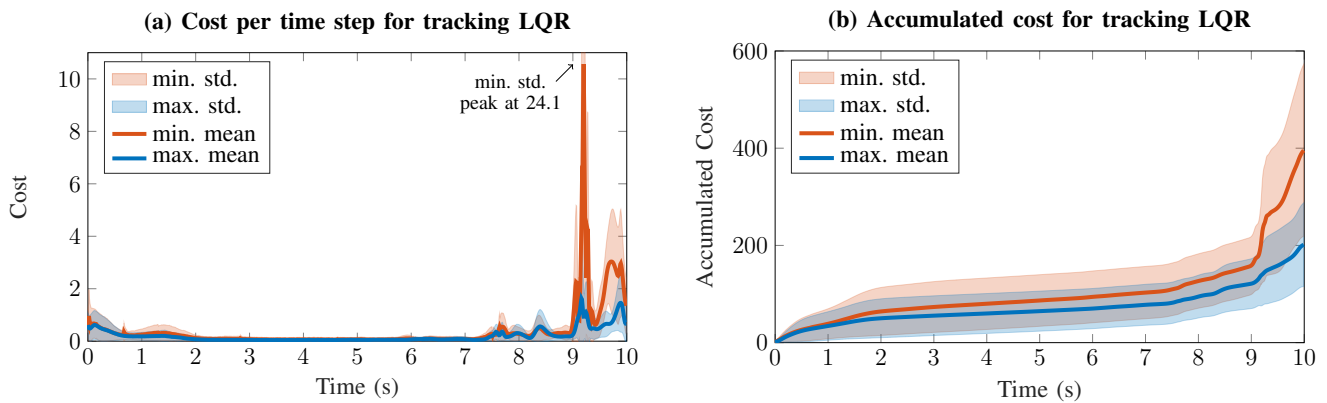
 
	\centering
	\begin{subfigure}{0.48\textwidth}
		\resizebox{\linewidth}{!}{\input{plots/triplecartpole.tex}}
		\refstepcounter{subfigure}\label{fig:cost}
	\end{subfigure}
	\quad
	\begin{subfigure}{0.48\textwidth}
		\resizebox{\linewidth}{!}{\input{plots/triplecartpolesum.tex}}
		\refstepcounter{subfigure}\label{fig:acost}
	\end{subfigure}
	\caption{Results for the tracking of the perturbed triple-pendulum cartpole swing-up with LQR. Maximal-coordinate mean cost (max. mean) and standard deviation (max. std.) in blue, minimal-coordinate mean cost (min. mean) and standard deviation (min. std.) in red. (a) Cost and standard deviation per time step. (b) Accumulated cost and corresponding standard deviation.}\label{fig:tracking}
\end{figure*} 

As disturbances, white Gaussian noise ($\mu=0,\,\sigma=2$) was added to the control input, viscous friction ($k=0.1$) to the prismatic and revolute joints, and uniformly distributed ($[0.9,1.1]$) perturbation factor was multiplied to the masses and inertias of cart and poles. Additionally, the initial configurations of cart and poles were uniformly distributed ($[-0.1,0.1]$) around the nominal initial configurations. For both controllers, 1000 runs were observed.

\subsection{Tracking LQR Results}

The results for the triple-pendulum cartpole tracking experiment are presented as cost per time step in Fig. \ref{fig:cost} and accumulated cost in Fig. \ref{fig:acost}. Overall, the tracking of the perturbed system with minimal- and maximal-coordinate LQR shows some similarities. For the first 7 seconds, energy is supplied to the system by swinging back and forth. During this period, both controllers manage to keep the incurred cost low. After 7 seconds, the final swing-up into the upright position takes place. An overall increase in cost per time step during the second phase becomes visible in Fig. \ref{fig:cost}, and especially the minimal-coordinate controller incurs increased cost. Additionally, a significantly higher cost variance for minimal-coordinate tracking can be observed in the final swing-up phase, when both systems diverge further from the nominal trajectory. Similarly, the accumulated costs for minimal- and maximal-coordinate LQR start to diverge visibly during the final swing-up phase in Fig. \ref{fig:acost}. Additionally, for the entire run the cost variance of minimal-coordinate tracking is larger than the maximal-coordinate cost variance.

\section{Discussion}\label{sec:discussion}
The investigation of the resulting control laws for the simple pendulum in Sec. \ref{sec:casestudy} demonstrated that there is a difference between minimal- and maximal-coordinate LQR control laws, even for matching cost functions, and that the maximal-coordinate control law can be interpreted as nonlinear, enabling potentially better performance. This understanding is supported by the simulation analysis, where, as expected, both minimal- and maximal-coordinate LQR are capable of stabilizing the systems under investigation from initial configurations close to the reference point or trajectory. However, for larger deviations from the nominal trajectory, maximal-coordinate LQR appears to be more robust. This interpretation stems from the larger basins of attraction for acrobot, cartpole, and the delta robot. And the tracking LQR problem also reveals interesting insights: During the first phase, both systems only deviated slightly from the reference trajectory and, therefore, both control laws appear to be reasonably accurate. In the second phase, with overall larger deviations from the nominal trajectory, the maximal-coordinate LQR performed better, leading to overall lower cost and cost variance. 

It is important to note, however, that further analysis is necessary to generalize our results. For certain nonlinear systems or certain regions, minimal-coordinate LQR could outperform maximal-coordinate LQR, for example, as seen with the pendulum in Fig. \ref{fig:respenderror} at certain velocities.

\section{Conclusions}\label{sec:conclusions}
We presented an extension of the linear-quadratic regulator to systems with physical constraints represented in maximal coordinates. Our derivation for such systems retains linearized constraints explicitly and directly incorporates them into the feedback gain calculations. Since the LQR derivation is built on the standard Riccati recursion, extensions such as state or input constraints are applicable to our derivation as well. The LQR controller can be directly deployed to stabilize systems around nominal trajectories.

The simulation of nonlinear systems suggests superior performance of maximal-coordinate LQR compared to minimal-coordinate LQR for comparable cost functions. This conclusion is based on the larger basins of attraction of maximal-coordinate LQR for the analyzed systems which generally points at a more robust controller. Time-varying tracking LQR shows comparable behavior for minimal- and maximal-coordinate LQR when close to the nominal trajectory. However, maximal-coordinate control appears to be less sensitive to larger deviations. Analyzing the control laws for a simple pendulum shows that the nonlinearity of the maximal-coordinate control law in minimal coordinates is one possible explanation for these results.

Overall, the results make a strong case for using maximal coordinates for linear-quadratic control, even for unconstrained systems. The ease of incorporating additional physical constraints during the modeling process adds to this conclusion. An interesting direction for future work could be to exploit sparsity when deriving the constrained LQR controller with additional Lagrange multipliers instead of, or in addition to, applying the more direct approach proposed in this paper. The maximal-coordinate LQR approach could also be extended to trajectory optimization methods such as iterative LQR (iLQR) or model predictive control (MPC) schemes based on LQR. 

\IEEEtriggeratref{17}
\bibliography{bib}

\begin{thebibliography}{10}

\bibitem{baraff_linear-time_1996}
D.~Baraff, ``Linear-time dynamics using {Lagrange} multipliers,'' in {\em
  Conference on {Computer} graphics and interactive techniques - {SIGGRAPH}
  '96}, pp.~137--146, ACM Press, 1996.

\bibitem{brudigam_linear-time_2020}
J.~Br\"udigam and Z.~Manchester, ``Linear-{Time} {Variational} {Integrators} in
  {Maximal} {Coordinates},'' in {\em Workshop on the {Algorithmic}
  {Foundations} of {Robotics} ({WAFR})}, Springer, 2020.

\bibitem{abraham_model-based_2017}
I.~Abraham, G.~de~la Torre, and T.~Murphey, ``Model-{Based} {Control} {Using}
  {Koopman} {Operators},'' in {\em Robotics: {Science} and {Systems} {XIII}},
  Robotics: Science and Systems Foundation, 2017.

\bibitem{suh_surprising_2020}
H.~Suh and R.~Tedrake, ``The {Surprising} {Effectiveness} of {Linear} {Models}
  for {Visual} {Foresight} in {Object} {Pile} {Manipulation},'' in {\em
  Workshop on the {Algorithmic} {Foundations} of {Robotics} ({WAFR})},
  Springer, 2020.

\bibitem{bevanda_koopman_2021}
P.~Bevanda, S.~Sosnowski, and S.~Hirche, ``Koopman {Operator} {Dynamical}
  {Models}: {Learning}, {Analysis} and {Control},'' {\em arXiv e-prints},
  vol.~arXiv:2102.02522 [eess.SY], 2021.

\bibitem{mills_force_1989}
J.~Mills and A.~Goldenberg, ``Force and position control of manipulators during
  constrained motion tasks,'' {\em IEEE Transactions on Robotics and
  Automation}, vol.~5, no.~1, pp.~30--46, 1989.

\bibitem{khayati_multistage_2006}
K.~Khayati, P.~Bigras, and L.-A. Dessaint, ``A {Multistage} {Position}/{Force}
  {Control} for {Constrained} {Robotic} {Systems} {With} {Friction}:
  {Joint}-{Space} {Decomposition}, {Linearization}, and {Multiobjective}
  {Observer}/{Controller} {Synthesis} {Using} {LMI} {Formalism},'' {\em IEEE
  Transactions on Industrial Electronics}, vol.~53, no.~5, pp.~1698--1712,
  2006.

\bibitem{townsend_optimal_1977}
M.~Townsend and T.~Tsai, ``On {Optimal} {Control} {Laws} for a {Class} of
  {Constrained} {Dynamical} {Systems} ({With} {Application} to {Control} of
  {Bipedal} {Locomotion}),'' {\em J. Dyn. Sys., Meas., Control}, vol.~99,
  no.~2, pp.~98--102, 1977.

\bibitem{posa_optimization_2016}
M.~Posa, S.~Kuindersma, and R.~Tedrake, ``Optimization and stabilization of
  trajectories for constrained dynamical systems,'' in {\em {International}
  {Conference} on {Robotics} and {Automation} ({ICRA})}, pp.~1366--1373, IEEE,
  2016.

\bibitem{you_tracking_1993}
L.-S. You and B.-S. Chen, ``Tracking control designs for both holonomic and
  non-holonomic constrained mechanical systems: a unified viewpoint,'' {\em
  Int. J. Control}, vol.~58, no.~3, pp.~587--612, 1993.

\bibitem{guanfeng_liu_unified_2002}
G.~Liu and Z.~Li, ``A unified geometric approach to modeling and control of
  constrained mechanical systems,'' {\em IEEE Transactions on Robotics and
  Automation}, vol.~18, no.~4, pp.~574--587, 2002.

\bibitem{kumar_feedback_1996}
A.~Kumar and P.~Daoutidis, ``Feedback regularization and control of nonlinear
  differential-algebraic-equation systems,'' {\em AIChE Journal}, vol.~42,
  no.~8, pp.~2175--2198, 1996.

\bibitem{krishnan_tracking_1994}
H.~Krishnan and N.~Mcclamroch, ``Tracking in nonlinear differential-algebraic
  control systems with applications to constrained robot systems,'' {\em
  Automatica}, vol.~30, no.~12, pp.~1885--1897, 1994.

\bibitem{kunkel_optimal_2008}
P.~Kunkel and V.~Mehrmann, ``Optimal control for unstructured nonlinear
  differential-algebraic equations of arbitrary index,'' {\em Math. Control.
  Signals, Syst.}, vol.~20, no.~3, pp.~227--269, 2008.

\bibitem{gerdts_direct_2003}
M.~Gerdts, ``Direct {Shooting} {Method} for the {Numerical} {Solution} of
  {Higher}-{Index} {DAE} {Optimal} {Control} {Problems},'' {\em Journal of
  Optimization Theory and Applications}, vol.~117, no.~2, pp.~267--294, 2003.

\bibitem{bender_linear-quadratic_1987}
D.~Bender and A.~Laub, ``The linear-quadratic optimal regulator for descriptor
  systems,'' {\em IEEE Transactions on Automatic Control}, vol.~32, no.~8,
  pp.~672--688, 1987.

\bibitem{kunkel_linear_1997}
P.~Kunkel and V.~Mehrmann, ``The linear quadratic optimal control problem for
  linear descriptor systems with variable coefficients,'' {\em Mathematics of
  Control, Signals, and Systems}, vol.~10, no.~3, pp.~247--264, 1997.

\bibitem{marz_linear_1995}
R.~M\"arz, ``On linear differential-algebraic equations and linearizations,''
  {\em Appl. Num. Mathematics}, vol.~18, no.~1-3, pp.~267--292, 1995.

\bibitem{kang_force_2003}
J.~Kang, S.~Bae, J.~Lee, and T.~Tak, ``Force {Equilibrium} {Approach} for
  {Linearization} of {Constrained} {Mechanical} {System} {Dynamics},'' {\em
  Journal of Mechanical Design}, vol.~125, no.~1, pp.~143--149, 2003.

\bibitem{negrut_practical_2006}
D.~Negrut and J.~Ortiz, ``A {Practical} {Approach} for the {Linearization} of
  the {Constrained} {Multibody} {Dynamics} {Equations},'' {\em J. Comput.
  Nonlinear Dynam.}, vol.~1, no.~3, pp.~230--239, 2006.

\bibitem{gonzalez_assessment_2017}
F.~Gonz\'alez, P.~Masarati, J.~Cuadrado, and M.~Naya, ``Assessment of
  {Linearization} {Approaches} for {Multibody} {Dynamics} {Formulations},''
  {\em J. Comput. Nonlinear Dynam.}, vol.~12, no.~4, p.~041009, 2017.

\bibitem{mason_balancing_2016}
S.~Mason, N.~Rotella, S.~Schaal, and L.~Righetti, ``Balancing and walking using
  full dynamics {LQR} control with contact constraints,'' in {\em
  {International} {Conference} on {Humanoid} {Robots} ({Humanoids})},
  pp.~63--68, IEEE, 2016.

\bibitem{savin_modification_2017}
S.~Savin, S.~Jatsun, and L.~Vorochaeva, ``Modification of constrained {LQR} for
  control of walking in-pipe robots,'' in {\em 2017 {Dynamics} of {Systems},
  {Mechanisms} and {Machines} ({Dynamics})}, pp.~1--6, IEEE, 2017.

\bibitem{johnson_structured_2015}
E.~Johnson, J.~Schultz, and T.~Murphey, ``Structured {Linearization} of
  {Discrete} {Mechanical} {Systems} for {Analysis} and {Optimal} {Control},''
  {\em IEEE Trans. Automat. Sci. Eng.}, vol.~12, no.~1, pp.~140--152, 2015.

\bibitem{scokaert_constrained_1998}
P.~Scokaert and J.~Rawlings, ``Constrained linear quadratic regulation,'' {\em
  IEEE Trans. Automat. Contr.}, vol.~43, no.~8, pp.~1163--1169, 1998.

\bibitem{mare_solution_2007}
J.~Mare and J.~De~Don\'a, ``Solution of the input-constrained {LQR} problem
  using dynamic programming,'' {\em Systems \& Control Letters}, vol.~56,
  no.~5, pp.~342--348, 2007.

\bibitem{nguyen_fast_2016}
H.-N. Nguyen and P.-O. Gutman, ``Fast {Constrained} {LQR} {Based} on {MPC}
  {With} {Linear} {Decomposition},'' {\em IEEE Trans. Automat. Contr.},
  vol.~61, no.~9, pp.~2585--2590, 2016.

\bibitem{ferranti_constrained_2016}
L.~Ferranti, G.~Stathopoulos, C.~Jones, and T.~Keviczky, ``Constrained {LQR}
  using online decomposition techniques,'' in {\em {Conference} on {Decision}
  and {Control} ({CDC})}, pp.~2339--2344, IEEE, 2016.

\bibitem{laine_efficient_2019}
F.~Laine and C.~Tomlin, ``Efficient {Computation} of {Feedback} {Control} for
  {Equality}-{Constrained} {LQR},'' in {\em {International} {Conference} on
  {Robotics} and {Automation} ({ICRA})}, pp.~6748--6754, IEEE, 2019.

\bibitem{knemeyer_minor_2020}
A.~Knemeyer, S.~Shield, and A.~Patel, ``Minor {Change}, {Major} {Gains}: {The}
  {Effect} of {Orientation} {Formulation} on {Solving} {Time} for
  {Multi}-{Body} {Trajectory} {Optimization},'' {\em IEEE Robotics and
  Automation Letters}, vol.~5, no.~4, pp.~5331--5338, 2020.

\bibitem{yilin_zhao_kinematics_1992}
Y.~Zhao and S.~BeMent, ``Kinematics, dynamics and control of wheeled mobile
  robots,'' in {\em {International} {Conference} on {Robotics} and {Automation}
  ({ICRA})}, pp.~91--96, IEEE Comput. Soc. Press, 1992.

\bibitem{mellinger_minimum_2011}
D.~Mellinger and V.~Kumar, ``Minimum snap trajectory generation and control for
  quadrotors,'' in {\em {International} {Conference} on {Robotics} and
  {Automation} ({ICRA})}, pp.~2520--2525, IEEE, 2011.

\bibitem{peng_constrained_2019}
Z.~Peng, J.~Wang, and J.~Wang, ``Constrained {Control} of {Autonomous}
  {Underwater} {Vehicles} {Based} on {Command} {Optimization} and {Disturbance}
  {Estimation},'' {\em IEEE Transactions on Industrial Electronics}, vol.~66,
  no.~5, pp.~3627--3635, 2019.

\bibitem{lewis_high_1997}
M.~Lewis and K.-H. Tan, ``High {Precision} {Formation} {Control} of {Mobile}
  {Robots} {Using} {Virtual} {Structures},'' {\em Autonomous Robots}, vol.~4,
  no.~4, pp.~387--403, 1997.

\bibitem{deshpande_distributed_2011}
P.~Deshpande, P.~Menon, C.~Edwards, and I.~Postlethwaite, ``A distributed
  control law with guaranteed {LQR} cost for identical dynamically coupled
  linear systems,'' in {\em Proceedings of the 2011 {American} {Control}
  {Conference} ({ACC})}, pp.~5342--5347, IEEE, 2011.

\bibitem{bertsekas_dynamic_1995}
D.~Bertsekas, {\em Dynamic programming and optimal control}.
\newblock Athena Scientific, 1995.

\bibitem{howell_altro_2019}
T.~Howell, B.~Jackson, and Z.~Manchester, ``{ALTRO}: {A} {Fast} {Solver} for
  {Constrained} {Trajectory} {Optimization},'' in {\em {International}
  {Conference} on {Intelligent} {Robots} and {Systems} ({IROS})},
  pp.~7674--7679, IEEE, 2019.

\end{thebibliography}
\bibliographystyle{ieeetr}

\end{document}